\documentstyle[11pt,amssymb,amstex]{amsart}

\def\bn{{\mathbb N}}

\def\l{\lambda} 

\def\m{\mu}

\def\w{\omega} \def\O{\Omega}

\newtheorem{thm}{Theorem}[section]

\newtheorem{cor}[thm]{Corollary}
\newtheorem{prop}[thm]{Proposition}

\theoremstyle{remark}

\def\id{{\bf 1}\!\!{\rm I}}
\def\nb{\nabla}

\begin{document}

\title[Weighted ergodic theorem ]{Weighted ergodic theorem for
contractions of Orlicz-Kantorovich lattice
$L_{M}(\widehat{\nabla},\widehat{\mu}$)}
\author{Inomjon Ganiev}
\address{Inomjon Ganiev\\
 Department of Science in Engineering\\
Faculty of Engineering, International Islamic University Malaysia\\
P.O. Box 10, 50728\\
Kuala-Lumpur, Malaysia}  \email{{\tt inam@@iium.edu.my},{\tt
ganiev1@@rambler.ru}}

\author{Farrukh Mukhamedov}
\address{Farrukh Mukhamedov\\
 Department of Computational \& Theoretical Sciences\\
Faculty of Science, International Islamic University Malaysia\\
P.O. Box, 141, 25710, Kuantan\\
Pahang, Malaysia} \email{{\tt farrukh\_m@@iium.edu.my}, {\tt
far75m@@yandex.ru}}

\begin{abstract}
In the present paper we prove Besocovich weighted ergodic theorem
for positive contractions acting on Orlich-Kantorovich space. Our
main tool is the use of methods of measurable bundles of
Banach-Kantorovich lattices.

\end{abstract}

\maketitle

\section{Introduction}

It is known that a pioneering work of von Neumman \cite{vN}
stimulated the development of the theory of Banach bundles (see
\cite{G2}). It was proved such a theory has vast applications in
analysis. Moreover, such a theory is well connected with
vector-valued Banach spaces, which has several application (see for
example, \cite{LL}). We recall that in the theory of Banach bundles
$L_0$-valued Banach spaces are considered, and such spaces are
called {\it Banach--Kantorovich spaces}. In \cite{G1,G2,K2}) the
theory of Banach--Kantorovich spaces were developed. Analogous of
many well known functional spaces has been defined and studied. For
example, in \cite{Ga3} Banach-Kantorovich lattice
$L_{p}(\widehat{\nabla},\widehat{\mu})$ is represented as a
measurable bundle of classical $L_{p}$ --lattices. In \cite{Z1,Z2}
an analogous of the Orlicz spaces has been considered. Naturally,
these functional Kantorovich spaces should have many of similar
properties like the classical ones, constructed by the real valued
measures.

We note that in \cite{BO} (see also \cite{Kr}) weighted ergodic
theorems for Dunford-Schwraz operators acting on $L_p$-spaces were
proved. Further, in \cite{B1,B2} such results were extended to
Banach-valued functions.
 In \cite{LW} it has been considered weighted ergodic
theorems and strong laws of large numbers. In \cite{W} some
properties of the convergence of Banach-valued martingales were
described and their connections with the geometrical properties of
Banach spaces were established too. Therefore, with the development
of the theory Banach--Kantorovich spaces
 there naturally arises the necessity to
study some ergodic type theorems for positive contractions and
martingales defined on such spaces.

To investigate the properties of Banach--Kantorovich spaces is
naturally to use measurable bundles of such spaces. Since, one has a
sufficiently well explored theory of measurable bundles of Banach
lattices \cite{G1}. Hence, it is an effective tool which gives well
opportunity to obtain various properties of Banach--Kantorovich
spaces \cite{Ga1},\cite{Ga2}. It is worth to mention that using this
way, in \cite{CGa,GaM2} weighted ergodic theorems for positive
contractions of Banach-Kantorovich lattices
$L_{p}(\widehat{\nabla},\widehat{\mu})$, has been established. In
\cite{Ga2} the convergence of martingales on such lattices is
proved. Further, in \cite{GaM} the "zero-two" law for positive
contractions of Banach-Kantorovich lattice
$L_{p}(\widehat{\nabla},\widehat{\mu})$ has been proved.

In \cite{ZC} individual ergodic theorem has been proved for positive
contractions of Orlicz-Kantorovich lattices.  In the present paper
we are going to prove Besicovich weighted ergodic theorem for
positive contractions acting on Orlich-Kantorovich space. Our
results extend and improve the results of \cite{GaM2,ZC}. To prove
the main result of this paper we are going to use measurable bundles
of Banach--Kantorovich lattices. We note that more effective methods
to study of Banach-Kantorovich spaces are the methods of
Boolean-valued analysis and measurable bundles (see \cite{K2}).

\section{Preliminaries}

In this section we recall necessary definitions and results of the
Banach-Kantorovich lattices.

 Let $(\Omega,\Sigma,\lambda)$ be a
measurable space with finite measure $\lambda$, and  $L_0(\Omega)$
be the algebra of all measurable functions on $\O$ ( here the
functions equal a.e. are identified) and let $\nb(\O)$ be the
Boolean algebra of all idempotents in $L_0(\Omega)$. By $\nb$ we
denote an arbitrary complete Boolean subalgebra of $\nb(\O)$.  By
${\cal L^{\infty}}(\Omega)$ we denote the set of all measurable
essentially bounded functions on $\O$, and $L^{\infty}(\Omega)$
denote an algebra of equivalence classes of essentially bounded
measurable functions.

Let $E$ be a linear space over the real field $\mathbb{R}$. By
$\|\cdot\|$ we denote a $L_0(\Omega)$-valued norm on $E$. Then the
pair $(E,\|\cdot\|)$ is called a {\it lattice-normed space (LNS)
over $L_0(\Omega)$}. An LNS $E$ is said to be {\it $d$-decomposable}
if for every $x\in E$ and the decomposition $\|x\|=f+g$ with $f$ and
$g$ disjoint positive elements in $L_0(\Omega)$ there exist $y,z\in
E$ such that $x=y+z$ with $\|y\|=f$, $\|z\|=g$.

Suppose that $(E,\|\cdot\|)$ is an LNS over $L_0(\O)$. A net
$\{x_\alpha\}$ of elements of $E$ is said to be {\it
$(bo)$-converging} to $x\in E$ (in this case we write $x=(bo)$-$\lim
x_\alpha$), if the net $\{\|x_\alpha - x\|\}$ $(o)$-converges to
zero (here $(o)$-convergence means the order convergence) in
$L_0(\Omega)$ (written as $(o)$-$\lim \|x_\alpha -x\|=0$). A net
$\{x_\alpha\}_{\alpha\in A}$ is called {\it $(bo)$-fundamental} if
$(x_\alpha-x_\beta)_{(\alpha,\beta)\in A\times A}$ $(bo)$-converges
to zero.

An LNS in which every $(bo)$-fundamental net $(bo)$-converges is
called {\it $(bo)$-complete}. A {\it Banach-Kantorovich space (BKS)
over $L_0(\Omega)$} is a $(bo)$-complete $d$-decomposable LNS over
$L_0(\Omega)$. It is well known \cite{K1,K2} that every BKS $E$ over
$L_0(\Omega)$ admits an $L_0(\Omega)$-module structure such that
$\|fx\|=|f|\cdot\|x\|$ for every $x\in E,\ f\in L_0(\Omega)$, where
$|f|$ is the modulus of a function $f\in L_0(\Omega)$.
 A BKS $({\cal U},\|\cdot\|)$ is called a {\it
Banach-Kantorovich lattice} if  ${\cal U}$ is a vector lattice and
the norm $\|\cdot\|$ is monotone, i.e.  $|u_1|\leq|u_2|$ implies
$\|u_1\|\leq\|u_2\|$. It is known \cite{K1} that the cone ${\cal
U}_+$ of  positive elements is $(bo)$-closed.

Let $(\O,\Sigma,\lambda)$ be the same as above and  $X$ be an
assisting real Banach space $(X(\omega),\|\cdot\|_{X(\omega)})$ to
each point $\omega\in\Omega$, where $ X(\omega)\neq\{0\}$ for all
$\w\in\O$. A {\it section} of $X$ is a function $u$ defined
$\l$-almost everywhere in $\O$ that takes values $u(\w)\in X(\w)$
for all $\w$ in the domain $dom(u)$ of $u$. Let $L$ be a set of
sections. The pair $(X,L)$ is called a {\it measurable Banach bundle
over $\Omega$} if
\begin{enumerate}
   \item[(1)] $\alpha_1 u_1 +\alpha_2 u_2 \in L$ for every $\alpha_1,\alpha_2\in\mathbb{R}$
   and $u_1,u_2\in L$, where\\
 $\alpha_1 u_1 + \alpha_2 u_2 : \omega\in dom(u_1) \cap
dom(u_2) \to \alpha_1 u_1(\omega) + \alpha_2 u_2(\omega)$;
   \item[(2)] the function $\|u\| : \omega\in dom(u) \to
\|u(\omega)\|_{X(\omega)}$ is measurable for every $u\in L$;
    \item[(3)] the set $\{u(\omega) : u\in L, \omega\in dom(u)\}$ is dense in
    $X(\omega)$ for every $\omega\in\O$.
\end{enumerate}

A measurable Banach bundle $(X,L)$ is called {\it measurable bundle
of Banach lattices (MBBL)} if $(X(\omega),\|\cdot\|_{X(\omega)})$ is
a Banach lattice for all $\w\in\O$ and for every $u_1,u_2\in L$ one
has $u_1\vee u_2\in L$, where $u_1 \vee u_2$: $\omega\in{\rm dom}\
(u_1)\cap {\rm dom}\
 (u_2) \rightarrow u_1(\omega)\vee u_2(\omega)$.

A section $s$ is called {\it step-section} if it has a form
$$
s(\omega)=\sum\limits_{i=1}^n \chi_{A_i}(\omega) u_i(\omega),
$$
for some $u_i\in L$, $A_i\in\Sigma$, $A_i\cap A_j=\emptyset$,
$i\neq j$, $i,j=1,\cdots,n$, $n\in\mathbb{N}$, where $\chi_A$ is
the indicator of a set $A$. A section $u$ is called {\it
measurable}  there exists a sequence of step-functions $\{s_n\}$
such that $s_n(\omega)\to u(\w)$ $\l$-a.e.

By $M(\Omega,X)$ we denote the set all measurable sections, and by
$L_0(\Omega,X)$ the factor space of $M(\Omega,X)$ with respect to
the equivalence relation of the equality a.e. Clearly, $L_0(\O,X)$
is an $L_0(\O)$-module. The equivalence class of an element $u\in
M(\O,X)$ is denoted by  $\hat{u}$. The norm of $\hat{u} L_0(\O,X)$
is defined as a class of equivalence in $L_0(\O)$ containing the
function $\|u(\omega)\|_{X(\omega)}$, namely
$\|\hat{u}\|=\widehat{(\|u(\w)\|_{X(\w)})}$. In \cite{G1} it was
proved that $L_0(\Omega,X)$ is a BKS over $L_0(\O)$. Furthermore,
for every BKS $E$ over $L_0(\O)$ there exists a measurable Banach
bundle $(X,L)$ over $\O$ such that $E$ is isomorphic to
$L_0(\Omega,X)$.

Let $X$ be a MBBL. We put $\hat u\leq\hat v$ if $u(\w)\leq v(\w)$
a.e. One can see that the relation  $\hat u\leq\hat v$ is a partial
order in  $L_0(\Omega,X)$. If $X$ is a MBBL, then  $L_0(\Omega,X)$
is a Banach-Kantarovich lattice \cite{Ga1,Ga3}.

A mapping $\mu : {\nabla}\to L_0(\Omega)$ is called a {\it
$L_0(\Omega)$-valued measure} if the following conditions are
satisfied:
\begin{enumerate}
\item[1)] $\mu(e)\geq 0$ for all $e\in{\nabla}$;
\item[2)] if $\ e\wedge g=0, e,g\in{\nabla},$ then $\mu(e\vee g)=\mu(e)+\mu(g)$;
\item[3)] if $e_n\downarrow 0,\ e_n\in{\nabla},\ n\in\mathbb{N}$, then
$\mu(e_n)\downarrow 0$.
\end{enumerate}

A $L_0(\Omega)$-valued measure $\mu$ is called {\it strictly
positive} if $\allowbreak \mu(e)=\nobreak 0,\ \allowbreak
e\in{\nabla}$ implies $e=0$.

Let a Boolean algebra ${\nabla}(\Omega)$ of all idempotents of
$L_0(\Omega)$ is a regular subalgebra of
 ${\nabla}$.

In the sequel we will consider a strictly positive
$L_0(\Omega)$-valued measure $\m$ with the following property
$\m(ge)=g\m(e)$ for all $e\in\nb$ and $g\in\nb(\O)$.

Let $\nabla_\omega,\ \omega\in\Omega$ be complete Boolean algebras
with strictly positive real-valued measures $\mu_\omega$. Put
 $\rho_\omega(e,g)=\mu_\omega(e\vartriangle g)$,
 $e,g\in\nabla_\omega$. Then $(\nabla_\omega,\mu_\omega)$ is a
 complete metric space. Let us consider a mapping
 $\nabla$, which assigns to each $\omega\in\Omega$ a Boolean algebra
 $\nabla_\omega$. Such a mapping is called a section.

 Assume that $L$ is a nonempty set of sections $\nabla$.
A pair $(\nabla,L)$ is called {\it a measurable bundle of Boolean
algebras over $\O$} if one has
\begin{enumerate}
\item[1)] $(\nabla,L)$  is a measurable bundle of metric spaces (see  \cite{Ga3});
\item[2)] if $e\in L$, then $e^{\perp}\in L$, where $e^{\perp} : \omega\in{\rm dom}\ (e) \rightarrow e^{\perp}(\omega)$;
\item[3)] if $e_1,e_2\in L$, then $e_1\vee e_2\in L$, where
$e_1\vee e_2 : \omega\in{\rm dom}\ (e_1) \cap {\rm dom}\ (e_2)
\rightarrow e_1(\omega) \vee e_2(\omega)$.
\end{enumerate}

Let $M(\Omega,\nabla)$ be the set of all measurable sections, and
$\hat{\nabla}$ be a factorization of $M(\Omega,\nabla)$ with respect
to equivalence relation the equality a.e. Let us define a mapping
$\hat{\mu}:\hat{\nabla}\rightarrow L_0(\Omega)$ by
$\hat{\mu}(\hat{e})=\hat{f}$, where $\hat{f}$ is a class containing
the function $f(\omega)=\mu_\omega(e(\omega))$. It is clear that the
mapping $\hat{\mu}$ is well-defined. It is known that
$(\hat{\nabla},\hat{\mu})$ is a complete Boolean algebra with a
strictly positive $L_0(\Omega)$-valued measure $\hat{\mu}$. Note
that a Boolean algebra ${\nabla}(\Omega)$ of all idempotents of
$L_0(\Omega)$ is identified with a regular subalgebra of
 $\hat{\nabla}$, and one has $\widehat{\mu}(g\widehat{e})=g\widehat{\mu}(\widehat{e})$ for all
 $g\in\nabla(\Omega)$ and $\widehat{e}\in\widehat{\nabla}$.

 The reverse is also tree, namely one has the following

\begin{thm}\cite{Ga3} Let $\tilde{\nabla}$ be a complete Boolean algebra, $\tilde{\mu}$ be a strictly positive
$L_0(\Omega)$-valued measure on $\tilde{\nabla}$, and
$\nabla(\Omega)$ is a regular subalgebra of $\tilde{\nabla}$ and
 $\tilde{\mu}(g\tilde{e})=g\tilde{\mu}(\tilde{e})$ for all $g\in\nabla(\Omega)$, $\tilde{e}\in\tilde{\nabla}$.
 Then there exists a measurable bundle of Boolean algebras
 $(\nabla,L)$ such that $\hat{\nabla}$ is isometrically isomorphic to
 $\tilde{\nabla}$.
 \end{thm}

By $L_0(\hat{\nabla},\hat{\mu})$ we denote an order complete vector
lattice $C_{\infty}(Q(\hat{\nabla})),$ where $Q(\hat{\nabla})$ is
the Stonian compact associated with complete Boolean algebra
$\hat{\nabla}$.  For $\widehat{f},\widehat{g}\in
L_0(\hat{\nabla},\hat{\mu})$ we let
 $\widehat{\rho}(\widehat{f},\widehat{g})=\allowbreak
 \int \frac{|\hat{f}-\hat{g}|}{1+|\hat{f}-\hat{g}|}d\hat{\mu}.$
 Then it is know \cite{Ga3} that $\widehat{\rho}$ is an
 $L_0(\Omega)$-valued metric on  $L_0(\hat{\nabla},\hat{\mu})$ and
 $(L_0(\hat{\nabla},\hat{\mu}),\widehat{\rho})$ is isometrically isomorphic to the measurable bundle
 of metric spaces $L_0(\nabla_\omega,\mu_\omega),$ where $\rho_\omega(a,b)= \int
 \frac{|a-b|}{1_{\omega}+|a-b|}d\mu_\omega.$ In particularly, each element
 $\widehat{f}\in L_0(\hat{\nabla},\hat{\mu})$ can be identified with the measurable section $\{f(\omega)\}_{\omega\in
 \Omega},$ here $f(\omega)\in L_0(\nabla_\omega,\mu_\omega).$

Following  the well known scheme of the construction of
$L_p$-spaces,  a space $L_p(\nb,\m)$ can be defined by
$$
L_p(\hat{\nb},\hat{\m})=\left\{\hat f\in L_0(\hat\nabla):
\int|\hat f|^pd\hat{\m} - \textrm{exist} \ \right\}, \ \ \ p\geq 1
$$
where $\hat\m$ is a  $L_0(\O)$-valued measure on $\hat\nb$.

It is known \cite{K1} that $L_p(\hat\nb,\hat\m)$ is a BKS over
$L_0(\O)$ with respect to the $L_0(\O)$-valued norm $\|\hat
f\|_{L_p(\hat{\nb},\hat{\m})}=\bigg(\int|\hat
f|^pd\hat{\m}\bigg)^{1/p}$. Moreover, $L_p(\hat{\nb},\hat{\m})$ is a
Banach-Kantorovich lattice (see \cite{K2},\cite{Ga3}).

Let $X$ be a mapping assisting an $L_p$-space constructed by a
real-valued measure $\m_\omega$, i.e.
$L_p(\nabla_\omega,\mu_\omega)$ to each point $\omega\in\O$ and let
$$
L=\bigg\{\sum\limits_{i=1}^n \alpha_i e_i : \alpha_i\in
{\mathbb{R}}, \ \ e_i\in M(\Omega,\nabla),\ i=\overline{1,n},\
n\in\mathbb{N}\bigg\}$$ be a set of sections. In \cite{Ga3, GaC}
it has been established that the pair $(X,L)$ is a measurable
bundle of Banach lattices and $L_0(\Omega,X)$ is modulo ordered
isomorphic to $L_p(\hat\nabla,\hat\mu)$.

  An even
continuous convex function $M : R \rightarrow [0,\infty)$ is called
an $N$-function, if $\lim\limits_{t\rightarrow0}\frac{M(t)}{t}=0$
  and $\lim\limits_{t\rightarrow\infty}\frac{M(t)}{t}=\infty$.
  Every $N$-function $M$ has the
form $M(t) = \int\limits_{0}^{|t|} p(s)ds,$ where $p(t)$ is a
nondecreasing function that is positive for $t > 0,$
right-continuous for $t\geq0,$ and such that $p(0) = 0$ and
$\lim\limits_{t\rightarrow\infty} p(t) =\infty$. Put $q(s) :=
\sup\{t : p(t)\leq s\},$ $s \geq 0.$ The function $N(t) :=
\int\limits_{0}^{|t|}q(s)ds $ is an N-function which is called the
complementary $N$-function to $M$ (see \cite{KR}).

The set $$L^0_M := L^0_ M(\hat{\nabla},\hat{\mu}) := \{x\in
L_0(\hat{\nabla}) : M(x)\in L_1(\hat{\nabla},\hat{\mu})\}$$ is
called the Orlicz $L_0$-class, and the vector space $$L_M :=
L_M(\hat{\nabla},\hat{\mu}) := \{x \in L_0(\hat{\nabla},\hat{\mu}) :
xy \in L_1(\hat{\nabla},\hat{\mu}) for\quad all\quad y\in L^0_N\}$$
is called the Orlicz $L_0$-space, where $N$ is the complementary
$N$-function to $M.$

 The following are valid: $L_M(\hat{\nabla},\hat{\mu}) \subset L_1(\hat{\nabla},\hat{\mu})$.

 Define the $L_0$-valued Orlicz norm on
$L_M(\hat{\nabla},\hat{\mu})$ as follows $$ \|x\|_M := \sup\{|\int
xy d\hat{\mu}| : y\in A(N)\},  x\in L_M(\hat{\nabla},\hat{\mu}),$$
where $A(N) = \{y \in L^0_ N : \int N(y)d\hat{\mu}\leq{\bf1}\}$ .
The pair $(L_M(\hat{\nabla},\hat{\mu}),\|\cdot\|_M) $ is a
Banach--Kantorovich lattice which is called the Orlicz--Kantorovich
lattice associated with the $L_0$-valued measure \cite{ZC}.

 As in the case of classical Orlicz spaces, along with the Orlicz norm  $ \|\cdot\|_{(M)}$ on $L_M(\hat{\nabla},\hat{\mu})$, we may
consider the $L_0$-valued Luxemburg norm $ \|x\|_{(M)} := \inf\{
\lambda\in L_0:\int M(\lambda^{-1}x)d\hat{\mu}\leq{\bf1},  \lambda
\quad is\quad an \quad invertible\quad positive\quad element \}$;
moreover, the pair  $(L_M(\hat{\nabla},\hat{\mu}), \|\cdot\|_{(M)})$
is also a Banach-Kantorovich lattice \cite{Z} (see also
\cite{Z1,Z2}).

Let $(L_M(\hat{\nabla},\hat{\mu}), \|\cdot\|_{(M)})$ be a
Orlicz--Kantorovich lattice and
$T:L_M(\hat{\nabla},\hat{\mu}),\hat{\mu})\to
L_M(\hat{\nabla},\hat{\mu})$ be a linear mapping. As usual  we
will say that $T$ is {\it positive} if $T\hat{f}\geq 0$ whenever
$\hat{f}\geq 0$. We say that $T$ is a {\it $L_0(\Omega)$-bounded
mapping} if there exists a function $k\in L_0(\Omega)$ such that
$\|T\hat{f}\|_{L_M(\hat\nb,\hat\m)}\leq k
\|\hat{f}\|_{L_M(\hat\nb,\hat\m)}$ for all $\hat{f}\in
L_M(\hat\nabla,\hat\mu)$. For a such mapping we can define an
element of $L_0(\O)$ as follows
$$
\|T\| =\sup\limits_{\|\hat{f}\|_{L_M(\hat\nb,\hat\m)}\leq\id}
\|T\hat{f}\|_{L_M(\hat\nb,\hat\m)},
$$
which is called an {\it $L_0(\O)$-valued norm} of $T$. If
$\|T\hat{f}\|_{L_M(\hat\nb,\hat\m)}\leq
\|\hat{f}\|_{L_M(\hat\nb,\hat\m)}$ then a mapping $T$ is said to be
a {\it $L_M(\hat\nb,\hat\m)$ contraction}.

The set of all essentially bounded functions w.r.t. $\hat{f}$
taken from $L_0(\hat\nb,\hat\m)$ is denoted by
$L^\infty(\hat\nb,\hat\m)$.

\section{$(o)$-convergence}

In this section we provide some auxiliary facts related to
$(o)$-convergence of sequence $\hat{f}_n$ from
$L_0(\hat{\nabla},\hat{\mu})$ and $(o)$-convergence of the sequence
$\{f_n(\omega)\}$ from $L_0(\nabla_\omega,\mu_\omega)$.

\begin{thm}\label{2.3} Let $\widehat{f}_{n}\in
L_0(\hat{\nabla},\hat{\mu}).$ Then $\sup\limits_n
 \hat{f}_n$ exists in $L_0(\hat{\nabla},\hat{\mu})$ if and only
 if $\sup\limits_n
 {f}_n(\omega)$ exists in $L_0(\nabla_\omega,\mu_\omega)$ for
 a.e. $\omega\in\Omega$. In the later case, one has
 $(\sup\limits_n
 \hat{f}_n)(\omega)=\sup\limits_n
 {f}_n(\omega)$ for
 a.e. $\omega\in\Omega$.
\end{thm}

\begin{pf} Assume that  $g(\omega)=\sup\limits_n {f}_n(\omega)$ exists in $(L_0(\nabla_\omega,\mu_\omega)$ for
 a.e. $\omega\in\Omega$. Denote
 $\hat{g}_n=\sup\limits_{1\leq k\leq n}\hat{f}_k$ in $L_0(\hat{\nabla},\hat{\mu}).$
 Then ${g}_n(\omega)=\sup\limits_{1 \leq k\leq n}{f}_k(\omega)$ for
 a.e. $\omega\in\Omega$.

 Obviously, that ${g}_{n}(\omega)\uparrow{g}(\omega)$ as
 $n\rightarrow\infty$ for
 a.e. $\omega\in\Omega$. The relation ${g}_{n}(\omega)\uparrow{g}(\omega)$ implies that
 ${g}_{n}(\omega)\stackrel{\rho_\omega}{\rightarrow}
 {g}(\omega)$ for
 a.e. $\omega\in\Omega$, this means $g\in  M(\Omega,X)$ and
 $\hat{g}\in L_0(\hat{\nabla},\hat{\mu})$.

Let us prove that $\hat{g} =
 \sup\limits_{n} \hat{f}_{n}$ in  $L_0(\hat{\nabla},\hat{\mu})$. It is
 clear that ${g}(\omega)\geq {f}_{n}(\omega)$ for
 a.e. $\omega\in\Omega$. Therefore, $\widehat{{g}}\geq
 \widehat{{f}}_{n}$ for all $n\in
\bn$.

Let $\widehat{\varphi}\in L^0(\hat{\nabla},\hat{\mu})$ and
$\widehat{{\varphi}}\geq
 \widehat{{f}}_{n}$ for all ${n}\in
\bn.$ Then ${\varphi}(\omega)\geq {f}_{n}(\omega)$ for any ${n}\in
\bn$. Hence, ${\varphi}(\omega)\geq {g}(\omega),$ for a.e.
$\omega\in\Omega$, i.e. $\hat{\varphi}\geq
 \hat{g}$. This yields that $\hat{g} = \sup\limits_{n\in\bn}
 \hat{f}_{n}$.

Conversely, let us assume that there exists such $\hat{\psi}\in
 L_0(\hat{\nabla},\hat{\mu})$ that $\hat{\psi} = \sup\limits_{n\in\bn} \hat{f}_{n}
  = \sup\limits_{n\in\bn}
 \hat{g}_{n}$.

From  $g_{n}(\omega) =
 \sup\limits_{1\leq k\leq n} f_{k}(\omega)$ for  a.e.
 $\omega\in\Omega$, we find $\psi(\omega)\geq f_{n}(\omega)$ for all
 $n\in\bn$
 for  a.e. $\omega\in\Omega$.
 Hence, one gets $\psi(\omega)\geq \sup\limits_{n\in\bn} f_{n}(\omega)
 = \sup\limits_{n\in\bn} g_{n}(\omega)$ for  a.e. $\omega\in\Omega$ .
 As $\hat{g}_{n} \rightarrow\hat{\psi}$ in
 metric $\hat{\rho}$, then  $g_{n}(\omega)\rightarrow\psi(\omega)$ in
 metric $\rho_\omega$ for  a.e.
 $\omega\in\Omega$.

 Since $\{g_n(\omega)\}$ is incrasing then $\psi(\omega) =
 \sup\limits_{n\in\bn} g_{n}(\omega)$ for  a.e.
 $\omega\in\Omega$.
\end{pf}

 From this theorem immediately follows two corollaries.

\begin{cor}\label{2.4} Let
 $\{\hat{f}_n\}\subset L_0(\hat{\nabla},\hat{\mu})$. Then
$\inf\limits_{n\in\bn} \hat{f}_{n}$ exists in
$L_0(\hat{\nabla},\hat{\mu})$
 if and only
 if $\inf\limits_{n\in\bn}
 f(\omega)$ exists in $L_0(\nabla_\omega,\mu_\omega)$ for  a.e.
 $\omega\in\Omega$. In later case, one has $(\inf\limits_{n\in\bn}
 \hat{f}_{n})(\omega)= \inf\limits_{n\in\bn} f_{n}(\omega)$ for  a.e.
 $\omega\in\Omega$.
 \end{cor}

\begin{cor}\label{2.5} Let $\hat{f}_{n} \in
 L_0(\hat{\nabla},\hat{\mu})$. If $\hat{f}_{n} \stackrel{(o)}\rightarrow \hat{f}$
 for some $\hat{f}\in L_0(\hat{\nabla},\hat{\mu})$, then
 $f_{n}(\omega) \stackrel{(o)}\rightarrow f(\omega)$ in $L_0(\nabla_\omega,\mu_\omega)$
 for  a.e. $\omega\in\Omega$. Conversely, if $f_{n}(\omega)
  \stackrel{(o)}\rightarrow g(\omega)$ for some $g(\omega)\in
 L_0(\nabla_\omega,\mu_\omega)$ for  a.e. $\omega\in\Omega$,
 then $\hat{g}\in L_0(\hat{\nabla},\hat{\mu})$ and $\hat{f}_{n} \stackrel{(o)}\rightarrow
 \hat{g}$ in $L_0(\hat{\nabla},\hat{\mu})$.\\
\end{cor}

\section{Weighted ergodic theorem}

In this section we are going to prove the main result of the paper.
We will prove $(o)$-convergence of Besicovich weighted ergodic
averages for positive contractions of the Orlich-Kantorovich space.

\begin{prop}\label{TT} Let $M$ be an $N$-function, strictly
convex in some interval and $T:L_M(\hat{\nabla},\hat{\mu})\to
L_M(\hat{\nabla},\hat{\mu})$ be a positive $L_0$-linear operator
such that
\begin{itemize}

\item[(i)] $\int M(|T\hat{f}| d\hat{\mu}\leq\int M(|\hat{f}|)
d\hat{\mu}$;

\item[(ii)] $\|T\|_{L_1(\hat\nabla,\hat\mu)\rightarrow
L_1(\hat\nabla,\hat\mu)}\leq \bf1$;

\item[(iii)] there exists $h\in L_M(\hat{\nabla},\hat{\mu}), h\neq0,
h(\omega)\neq0 $ such that $Th=h$.
\end{itemize}
Then there exists a family $T_\omega :
L_M(\nabla_\omega,\mu_\omega)\to L_M(\nabla_\omega,\mu_\omega)$ of
positive linear operators such that for each $\hat{f}\in
L_M(\nabla_\omega,\mu_\omega)$ the equality $T_\omega f(\omega) =
(T\hat{f})(\omega)$  holds for almost all $\omega\in\Omega$.
Moreover, one has
\begin{itemize}

\item[(a)]  $\int M(|T_\omega f(\omega)| d{\mu}_\omega\leq\int
M(|f(\omega)|) d{\mu}_\omega$ for almost all $\omega\in\Omega$;

\item[(b)] $\|T_\omega\|_{L_M(\nabla_\omega,\mu_\omega)\rightarrow
L_M(\nabla_\omega,\mu_\omega)}\leq 1$;

\item[(c)] There exists $h(\omega)\in L_M(\nabla_\omega,\mu_\omega),
h(\omega)\neq0$ such that  $T_\omega h(\omega)= h(\omega)$.
\end{itemize}
\end{prop}

\begin{pf} Since $\|T\|_{L_1(\hat\nabla,\hat\mu)\rightarrow
L_1(\hat\nabla,\hat\mu)}\leq \bf1$ by Theorem 2.1 \cite{GaM} there
exists a family $T_\omega : L_1(\nabla_\omega,\mu_\omega)\to
L_1(\nabla_\omega,\mu_\omega)$ of positive linear operators such
that for each $\hat{f}\in L_M(\nabla_\omega,\mu_\omega)$ the
equality $T_\omega f(\omega) = (T\hat{f})(\omega)$  holds for almost
all $\omega\in\Omega$. Moreover, one has
$\|T_\omega\|_{L_1(\nabla_\omega,\mu_\omega)\rightarrow
L_1(\nabla_\omega,\mu_\omega)}\leq 1$.

(a) Let $\hat{f}=\{f(\omega)\}_{\omega\in\Omega}\in
L^\infty(\hat\nabla,\hat\mu).$ Then
\begin{eqnarray*}
\int M(|T_\omega f(\omega)|) d{\mu}_\omega&=&\int
M(|T\hat{f}|(\omega))  d{\mu_\omega}=\bigg(\int M(|T\hat{f}|
d\hat{\mu}\bigg)(\omega)\\[2mm]
&\leq&\bigg(\int M(|\hat{f}|) d\hat{\mu}\bigg)(\omega)=\int
M(|f(\omega)|) d{\mu}_\omega.
\end{eqnarray*}

Let $f(\omega)\in L_M(\nabla_\omega,\mu_\omega)$, then we choose a
sequence  $\hat{f}_n\in L^\infty(\hat\nabla,\hat\mu)$ such that
$f(\omega)=\lim\limits_{n\rightarrow\infty}f_n(\omega)$. Then
$T_\omega f(\omega)=\lim\limits_{n\rightarrow\infty}T_\omega
f_n(\omega)$ and
\begin{eqnarray*}
\int M(|T_\omega f(\omega)|) d{\mu}_\omega&=&
\lim\limits_{n\rightarrow\infty}\int M(|T_\omega f_n(\omega)|)
d{\mu}_\omega\\[2mm]
&=&\lim\limits_{n\rightarrow\infty}\bigg(\int M(|T\hat{f_n}|
d\hat{\mu}\bigg)(\omega)\\[2mm]
&\leq&\lim\limits_{n\rightarrow\infty}\bigg(\int M(|\hat{f_n}|
d\hat{\mu}\bigg)(\omega)\\[2mm]
&=&\lim\limits_{n\rightarrow\infty}\bigg(\int M(|f_n(\omega)|\bigg)
d\mu_\omega\\
&=&\int M(|f(\omega)|) d{\mu}_\omega.
\end{eqnarray*}

(c) Since, there exists $h\in L_M(\hat{\nabla},\hat{\mu}), h\neq0,
h(\omega)\neq0 $ such that $Th=h$ we have that $T_\omega
h(\omega)=(Th)(\omega)= h(\omega).$

(b)  Using (a), (c),
$\|T_\omega\|_{L_1(\nabla_\omega,\mu_\omega)\rightarrow
L_1(\nabla_\omega,\mu_\omega)}\leq 1$ and Theorem 2.5 \cite{Gal} we
have $\|T_\omega f(\omega)\|_\infty\leq \|f(\omega)\|_\infty $ for
every $f(\omega)\in L_1(\nabla_\omega,\mu_\omega)\cap
L^\infty(\nabla_\omega,\mu_\omega)$, consequently the last
inequality holds for every $f(\omega)\in
L_M(\nabla_\omega,\mu_\omega)\cap
L^\infty(\nabla_\omega,\mu_\omega)$. The positivity $T_\omega$
implies $T_\omega{\bf1}_\omega\leq {\bf1}_\omega$. Since
$\|T_\omega\|_{L_1(\nabla_\omega,\mu_\omega)\rightarrow
L_1(\nabla_\omega,\mu_\omega)}\leq 1$, $T_\omega{\bf1}_\omega\leq
{\bf1}_\omega$ and the Orlicz space $L_M(\nabla_\omega,\mu_\omega)$
is an interpolation space with interpolation constant unity (see
\cite{kr}, II, sec.4, item 6) we have
$T_\omega(L_M(\nabla_\omega,\mu_\omega))\subset
L_M(\nabla_\omega,\mu_\omega)$. This completes the proof.
\end{pf}

 A sequence $\{b(k)\}$ is called {\it Besicovich} if for
every $\varepsilon>0$ there is a sequence of trigonometric
polynomials $\psi_{\varepsilon}$, such that
$$\lim\limits_{N\rightarrow \infty}\sup\frac{1}{N}\sum\limits_{k=1}^{N}|b(k)-\psi_{\varepsilon}(k)|<\varepsilon$$

We say that $\{b(k)\}$ is {\it bounded Besicovich} if $b(k)\in
\ell^{\infty}.$ In the present paper we consider only bounded, real
Besicovich sequences.

\begin{thm}  Let $M$ be an $N$-function, strictly convex
in some interval and $T:L_M(\hat{\nabla},\hat{\mu})\to
L_M(\hat{\nabla},\hat{\mu})$ be a positive $L_0$-linear operator
such that
\begin{itemize}

\item[(i)] $\int M(|T\hat{f}| d\hat{\mu}\leq\int M(|\hat{f}|)
d\hat{\mu}$;

\item[(ii)] $\|T\|_{L_1(\hat\nabla,\hat\mu)\rightarrow
L_1(\hat\nabla,\hat\mu)}\leq \bf1$;

\item[(iii)] There exists $h\in L_M(\hat{\nabla},\hat{\mu}), h\neq0,
h(\omega)\neq0 $ such that $Th=h$.
\end{itemize}
Then
\begin{itemize}
\item[(a)] There exists
$$\sup\limits_{n\geq 1}\frac{1}{n}\sum\limits_{k=1}^{n-1}b(k)T^k|\hat{f}|$$
in $L_M(\hat{\nabla},\hat{\mu})$ for any bounded Besicovich sequence
$b(k)$.

\item[(b)] If $\{b(k)\}$ is a bounded  Besicovich sequence, then for
every $\hat{f}\in L_M(\hat{\nabla},\hat{\mu})$ there exists
$\hat{f}^*\in L_M(\hat{\nabla},\hat{\mu})$ such that the sequence
$$\widetilde{A_{n}}(\widehat{f})=\frac{1}{n}\sum\limits_{k=1}^{n-1}b(k)T^k\widehat{f}$$
$(o)$-converges in $L_M(\widehat{\nabla},\widehat{\mu}).$
\end{itemize}
\end{thm}

\begin{pf} (a). Let $\hat{f}\in L_M(\hat{\nabla},\hat{\mu})$
and $T:L_M(\hat{\nabla},\hat{\mu})\to L_M(\hat{\nabla},\hat{\mu})$
satisfy conditions (i),(ii),(iii). Then by Proposition \ref{TT}
there exists a family $T_\omega : L_M(\nabla_\omega,\mu_\omega)\to
L_M(\nabla_\omega,\mu_\omega)$ of positive linear operators such
that for each $\hat{f}\in L_M(\hat{\nabla},\hat{\mu})$ the equality
$T_\omega f(\omega) = (T\hat{f})(\omega)$  holds for almost all
$\omega\in\Omega$. Moreover one has
 $\int M(|T_\omega f(\omega)| d{\mu}_\omega\leq\int
M(|f(\omega)|) d{\mu}_\omega$ for almost all $\omega\in\Omega$, and
$\|T_\omega\|_{L_1(\nabla_\omega,\mu_\omega)\rightarrow
L_1(\nabla_\omega,\mu_\omega)}\leq 1$. Besides, there exists
$h(\omega)\in L_M(\nabla_\omega,\mu_\omega), h(\omega)\neq0$ such
that $T_\omega h(\omega)= h(\omega)$. This implies that the family
$\{T_\omega\}$ satisfies conditions of Theorem 3.1 \cite{Gal}.
Consequently, we get
\begin{eqnarray*}
\bigg\|\sup\limits_{n\geq
1}\frac{1}{n}\sum\limits_{k=1}^{n-1}b(k)T_\omega^k(|f(\omega)|)\bigg\|_{(M)}
&\leq&
 b \bigg\|\sup\limits_{n\geq 1}\frac{1}{n}\sum\limits_{k=1}^{n-1}T_\omega^k(|f(\omega)|)\bigg\|_{(M)}\\[3mm]
 &\leq& bC \|f(\omega)\|_{(M)},
\end{eqnarray*}
 where $b=\sup\limits_{k}|b_k|$.

 Therefore, $\sup\limits_{n\geq
 1}\frac{1}{n}\sum\limits_{k=1}^{n-1}b(k)T_\omega^k(|f(\omega)|)$
 exists in $L_M(\nabla_\omega,\mu_\omega).$ Thus, from Theorem
 4.1 \cite{Ga3} and Proposition 2.3 \cite{ZC} one finds
 the existence of
 $$\sup\limits_{n\geq 1}\frac{1}{n}\sum\limits_{k=1}^{n-1}b(k)T^k|\hat{f}|$$
in $L_M(\hat{\nabla},\hat{\mu})$.

(b). As the family $\{T_\omega\}$ satisfies conditions Theorem 3.1
\cite{Gal}, then there exists $f^*(\omega)\in
L_M(\nabla_\omega,\mu_\omega)$ such that
$$\lim\limits_{n\rightarrow\infty}{A_{n}}({f}(\omega))=\lim\limits_{n\rightarrow\infty}\frac{1}{n}\sum\limits_{k=1}^{n-1}b(k)T_\omega^k(f(\omega))=f^*(\omega)$$
$\mu_\omega$- almost everywhere  for each  $f(\omega)\in
L_M(\nabla_\omega,\mu_\omega)$. Taking into account that
$\mu_\omega$- almost everywhere convergence for sequence in
$L_0(\nabla_\omega,\mu_\omega)$ coincides with $(o)$-convergent in
$L_0(\nabla_\omega,\mu_\omega)$, we infer that the sequence
${A_{n}}({f}(\omega))$ is $(o)$-convergent in
$L_0(\nabla_\omega,\mu_\omega)$. Hence, by Corollary \ref{2.5} one
finds that $\widetilde{A_{n}}(\widehat{f})$ is $(o)$-convergent in
$L_0(\widehat{\nabla},\widehat{\mu})$.

The existence $\sup\limits_{n\geq
1}\frac{1}{n}\sum\limits_{k=1}^{n-1}b(k)T^k|\hat{f}|$ in
$L_M(\hat{\nabla},\hat{\mu})$ yields that
\begin{equation*}
\widetilde{A_{n}}(\widehat{f})=\widehat{{A_{n}}({f}(\omega)}){\stackrel{(o)}{\rightarrow}}\widehat{f^*}
=\widehat{{f^*}(\omega)}
\end{equation*}
in $L_M(\widehat{\nabla}, \widehat{\mu}).$ This completes the proof.
\end{pf}

\section*{Acknowledgement}
 The first  named author (I.G.) acknowledges the IIUM-RMC
grant  EDW B 11-185-0663. The second named author (F.M.)
acknowledges the MOHE Grant FRGS11-022-0170. He also thanks the
Junior Associate scheme of the Abdus Salam International Centre for
Theoretical Physics, Trieste, Italy.

\end{document}